\documentclass[12pt, letterpaper]{amsart}

\title{Gregory Trees, the Continuum, and Martin's Axiom}
\author{Kenneth Kunen}
\address{Department of Mathematics \\
University of Wisconsin \\
Madison, WI 57306 USA}
\email{kunen@math.wisc.edu}
\urladdr{http://www.math.wisc.edu/\symbol{126}kunen}
\author{Dilip Raghavan}
\thanks{Authors partially supported by NSF Grant DMS-0456653.}
\address{Department of Mathematics \\
University of Wisconsin \\
Madison, WI 57306 USA}
\email{raghavan@math.wisc.edu}
\urladdr{http://www.math.wisc.edu/\symbol{126}raghavan}

\date{\today}


\usepackage{amssymb, mathrsfs, latexsym, bbm}

\newtheorem{Theorem}{Theorem}[section]
\newtheorem{Lemma}[Theorem]{Lemma}
\newtheorem{Proposition}[Theorem]{Proposition}
\newtheorem{Def}[Theorem]{Definition}

\newtheorem{Question}[Theorem]{Question}

\newcommand\one{\mathbbm{1}}  
\newcommand\res{\mathord {\upharpoonright}}  
\newcommand{\PPP}{\mathbb{P}}
\newcommand{\QQQ}{\mathbb{Q}}
\newcommand{\BBB}{\mathbb{B}}
\newcommand{\RRR}{\mathbb{R}}
\newcommand{\SSS}{\mathbb{S}}
\renewcommand{\P}{\mathbb{P}}
\newcommand{\Q}{\mathbb{Q}}
\newcommand{\T}{\mathbb{T}}
\newcommand\Fn{\mathrm{Fn}} 
\newcommand\MA{\mathrm{MA}} 
\newcommand\PFA{\mathrm{PFA}} 
\newcommand\CH{\mathrm{CH}} 
\newcommand\cat{^{\mathord{\frown}}}  

\newcommand{\forces}{\Vdash}

\renewcommand{\c}{\mathfrak{c}}
\renewcommand{\[}{\left[}
\renewcommand{\]}{\right]}
\DeclareMathOperator{\cf}{cf}
\DeclareMathOperator{\hgt}{ht}

\newcommand{\V}{{\mathbf{V}}}
\newcommand{\LL}{{\mathbf{L}}}
\newcommand{\VG}{{{\mathbf{V}}[G]}}

\begin{document}
\begin{abstract}
We continue the investigation of Gregory trees and the Cantor Tree
Property carried out by Hart and Kunen.
We produce models of $\MA$ with the Continuum arbitrarily large
in which there are Gregory trees, and in
which there are no Gregory trees.
\end{abstract}
\maketitle
\section{Introduction}\label{sec:intro}

We view the tree $2^{<\omega_1}$ as a forcing poset,
defining $p \le q$ iff $p \supseteq q$; so
$\one = \emptyset$, the empty sequence.
A \emph{Gregory tree} is a type of subtree of $2^{<\omega_1}$
which is ``almost countably closed''.
The notion is due to Gregory \cite{Gr}, although the terminology
in the next definition is from Hart and Kunen \cite{KH}.

\begin{Def} \label{def:Gtrees}
A \emph{Cantor tree} in $2^{<\omega_1}$ is a subset
$\{f_{\sigma}: \sigma \in 2^{< \omega}\} \subseteq 2^{< \omega_1}$
such that for all $\sigma \in 2^{< \omega}$, $f_{{\sigma}^{\frown}0}$ and
$f_{{\sigma}^{\frown}1}$ are incompatible nodes in $2^{< \omega_1}$
that extend $f_{\sigma}$. A subtree $\T$ of $2^{<{\omega_1}}$ has
the \emph{Cantor Tree Property (CTP)} iff
\begin{itemize}
\item[1.]
For every $f \in \T$, $f^{\frown}0, f^{\frown}1 \in \T$.
\item[2.]
Given any Cantor tree $\{f_{\sigma}: \sigma \in 2^{< \omega}\} \subset \T$,
there are $x \in {2}^{\omega}$ and $g \in \T$ such that
$\forall n \in \omega \[ g \le {f}_{x \res n}  \] $. 
\end{itemize}
A subtree $\T$ of $2^{< \omega_1}$ is a \emph{Gregory tree} iff it has the CTP,
but does not have a cofinal branch.
\end{Def}

Paper \cite{KH2} relates Gregory trees to more general forcing posets
with the CTP.

\begin{Theorem}[Gregory \cite{Gr}] 
\label{thm:greg-tree}
$2^{\aleph_0} < 2^{\aleph_1}$ implies that there is a Gregory tree.
\end{Theorem}

Gregory trees are of interest in the theory of proper forcing.
It is easy to see (Lemma 5.5 of \cite{KH}) that a Gregory tree
$\T$ is a totally proper poset, that is, it is proper and does not
add any reals.  Moreover, forcing with $\T$ adds a cofinal branch through $\T$.
One might hope to do a countable support iteration of
these  totally proper forcings, producing a model of $\CH$ plus no
Gregory trees, but this is impossible by Theorem
\ref{thm:greg-tree}, so that the iteration must add reals,
although the CTP is annoyingly close to being countably closed.
Of course, 

\begin{Proposition}
\label{prop:no-greg-tree}
$\PFA$ implies that there are no Gregory trees.
\end{Proposition}

Gregory trees arose naturally in \cite{KH} from the study of
\emph{weird} topological spaces; these are compact non-scattered
(Hausdorff) spaces $X$ such that no perfect subset
$P \subseteq X$ is totally disconnected.
Assuming $\diamondsuit$, there is a weird space
which is hereditarily Lindel{\"o}f (HL); this is false under PFA,
since by Lemma 5.7 of \cite{KH},
if $X$ is compact, HL, and not totally disconnected,
and $X$ has no subspaces homeomorphic to the Cantor set,
then there is a Gregory tree. 

Whenever a result,
is proved from PFA, two natural questions arise. 
First, does it follow just from $\MA + \neg \CH$?
Second, is it consistent with $2^{\aleph_0} > \aleph_2$?
Of course, the second question is trivial if the answer to the
first question is ``yes''. In this paper, with regard to
Proposition \ref{prop:no-greg-tree},
we show that the answer is ``no'' to the first question and ``yes''
to the second.
In Section \ref{sec:MAplusnotCHplusexistsGtree},
we produce models of $\MA$ in which there exists a Gregory tree;
$\c$ can be ``anything regular''.
In Section, \ref{sec:noGtreesplusclarge}
we produce models of $\MA + \neg \CH$ in
which there does not exist a Gregory tree;
here, $\c$ can be ``anything regular'' \emph{except} the successor
of a cardinal of cofinality $\omega$, so we are left with the
following:
\begin{Question}
\label{ques:greg}
Assume $\MA$ and $\c = {\aleph}_{\omega + 1}$.
Must there be a Gregory tree?
\end{Question}

\section{Notational Conventions for Iterated c.c.c.\@ Forcing}
\label{sec:cccforcing}
In this paper we only consider finite support iterations of c.c.c.\@ forcings. Before giving the proofs of our theorems, we set out some notational conventions regarding these iterations.

As usual in forcing, a \textit{forcing poset} $\PPP$ really denotes a
triple, $(\PPP, \le, \one)$, where $\le$ is a transitive reflexive relation on $\PPP$ and $\one $ is a largest element of $\PPP$.  Then, the notation $\PPP \subseteq \QQQ$ implies that the orders agree and that $\one_\PPP = \one_\QQQ$. $\PPP \subseteq_c \QQQ$ means that in addition, $\PPP$ is 
a \emph{complete} sub-order of $\QQQ$; this implies that we may view the $\QQQ$--extension as a generic extension  of the $\PPP$--extension 
(see, e.g., \cite{Kunen}). Since all our iterated forcings are c.c.c.\@ with finite supports, it is simpler \emph{not} to follow the approach of \cite{Kunen}, but rather to construct in the ground model a \emph{normal chain} of c.c.c.\@ posets, $\langle \PPP_\alpha : \alpha \le \kappa \rangle$, where $\alpha < \beta \to \PPP_\alpha \subseteq_c \PPP_\beta$ and we take 
unions at limits (which preserves the c.c.c.).
In standard iterated forcing constructions, the $\P_\alpha$
are constructed inductively; given $\P_\alpha$, we choose
$\mathring{\QQQ}_\alpha$, which is a $\PPP_\alpha$--name 
forced by $\one$ to be a c.c.c.\@ poset;
then $\PPP_{\alpha+1}$ is identified with
$\PPP_\alpha *  \mathring{\QQQ}_\alpha$.
However, the basic theory of these iterations does not require a
$\mathring{\QQQ}_\alpha$; in Section \ref{sec:noGtreesplusclarge}, 
it will sometimes be convenient to view a $\gamma$--chain as
a $\cf(\gamma)$--chain by restricting to a cofinal sequence.

We shall always take $\PPP_0 = \{\one\}$, so that we can identify the $\PPP_0$--extension with the ground model. If $G$ is a $(\V, \PPP_\kappa)$--generic
filter, then $G_\alpha := G \cap \PPP_\alpha $ is $(\V, \PPP_\alpha)$--generic,
and we let ${\V}_\alpha = \V[G_\alpha]$; so, ${\V}_0 = \V$.

If $\varphi$ is a sentence in the $\PPP_\alpha$ forcing language and $p \in \PPP_\alpha$,
then $p \forces_\alpha \varphi$ abbreviates $p \forces_{\PPP_\alpha} \varphi$.  Note that we need a subscript on the $\forces$, since for any $\beta > \alpha$, the assertion ``$p \forces_\beta \varphi$'' is meaningful, although its truth can vary with $\beta$.	  

We use Shoenfield--style names as in \cite{Kunen}; that is, a name is
a set of ordered pairs of names and forcing conditions.
So, an inclusion of names ($\mathring A \subseteq \mathring B$)
implies an inclusion of the sets named
($\one \forces \mathring A \subseteq \mathring B$).
Also, if $\PPP_\alpha \subseteq_c \PPP_\beta$,
then every $\PPP_\alpha$--name is a $\PPP_\beta$--name.
In Section \ref{sec:MAplusnotCHplusexistsGtree},
we shall build a Gregory tree $\T$ in $\VG$ by constructing in $\V$
an ascending sequence of names
$\langle \mathring \T_\alpha : \alpha \le \kappa\rangle$,
where $\mathring \T_\alpha$ is a $\P_\alpha$--name.

If $G$ is $\P$--generic over $\V$ and $X \in \V$, then 
every subset of $X$ in $\VG$ is named by a \emph{nice name
$\mathring b$ for a subset of $X$}; so
$\mathring b = \bigcup\{ \{\check x\} \times E_x : x \in X \}$, where
each $E_x$ is an antichain in $\P$ (see \cite{Kunen}).
Also, if $p \forces \mathring a \subseteq \check X$
then there is a nice name $\mathring b$ for a subset of $X$ such
that $p \forces \mathring a  = \mathring b$.
With iterated forcing, where $\P = \P_\gamma$ results from a
normal chain of c.c.c.\@ posets
$\langle \PPP_\alpha : \alpha \le \gamma \rangle$:
if $\cf(\gamma) \ge \omega_1$ and $X$ is countable, then,
since the antichains are also countable,
there is an $\alpha_0 < \gamma$ such that our
$\mathring b$ is also a nice $\P_\alpha$ name
whenever $\alpha_0 \le \alpha \le \gamma$.

\section{A Model of $\MA + \neg \CH +  \ \text{There is a Gregory Tree}$}
\label{sec:MAplusnotCHplusexistsGtree} 
 
\begin{Theorem}
\label{Thm:is-greg-tree}
Assume that in the ground model $\V$, $\kappa \ge \aleph_2$ and
$\kappa^{<\kappa} = \kappa$. Then there is a c.c.c.\@ forcing
extension $\V[G]$ satisfying $\MA + 2^{\aleph_0} = \kappa$ in which
there is a Gregory tree.
\end{Theorem}
\begin{proof}
The standard procedure for constructing a model of $\MA$
in which some consequence of PFA fails is to start with a
counter-example in $\V$ which is not destroyed by the c.c.c.\@ iteration.
However, every Gregory tree $\T$ in $\V$ is
destroyed immediately whenever a real is added, since that will cause the CTP to fail. Instead, our tree $\T$ will grow along with the
iterated forcing which produces our model. To do this, we inductively construct
the following, satisfying the listed conditions:
	\begin{enumerate}
		\item 
			$\langle \PPP_\alpha : \alpha \le \kappa\rangle$ is a normal chain of c.c.c.\@ posets.
		\item 
			$|\PPP_\alpha| < \kappa$ for all $\alpha < \kappa$.
		\item 
			$\PPP_{\alpha + 1} \cong \PPP_\alpha * \Fn(\omega,2) * \mathring{\QQQ}_\alpha$, where $\one \forces_\alpha$ ``$\mathring{\QQQ}_\alpha$ is c.c.c.''.
		\item 
			Each $\mathring \T_\alpha$ is a $\PPP_\alpha$--name,
for $\alpha \le \kappa$.
		\item 
			$\one \forces_\alpha$ ``$\mathring \T_\alpha$ is a subtree of $2^{< \omega_1}$ and $\forall f \in \mathring \T_\alpha \forall s \in 2^{<\omega} \,[f \cat s \in \mathring \T_\alpha]$''.
		\item 
			$\mathring \T_0$ is a name for the Cantor tree $2^{< \omega}$.
		\item 
			If $\alpha < \beta$ then $\mathring \T_\alpha \subseteq \mathring \T_\beta$, so that $\one \forces_\beta \mathring \T_\alpha \subseteq \mathring \T_\beta$.
		\item 
			If $\gamma$ is a limit, then $\mathring \T_{\gamma} = \bigcup_{\alpha < \gamma} \mathring \T_\gamma$.
		\item 
			$\mathring g_\alpha$ is a $\PPP_{\alpha+1}$--name and $\one \forces_{\alpha+1}$ ``$\mathring \T_{\alpha+1} = \mathring \T_\alpha \cup \{\mathring g_\alpha \cat s : s \in 2^{< \omega}\} $''
		\item 
			$\one \forces_\alpha$ ``$\mathring \T_\alpha$ has no uncountable chains''.
		\item 
			$\one \forces_{\alpha+1}$ ``$\mathring \T_\alpha$ is special''.
		\item 
			$\mathring c_\alpha$ is a $\PPP_{\alpha+1}$--name for the function from $\omega$ to $2$ added by the $\Fn(\omega,2)$ in item (3).
		\item 
			$\mathring K_\alpha$ and $\mathring k_{\alpha,s}$ are $\PPP_\alpha$--names whenever $\alpha < \omega_1$ and $s \in 2^{<\omega}$.
		\item 
			$\one \forces_\alpha$ ``$\mathring K_\alpha$ is a Cantor tree in $\mathring \T_\alpha$, and $\mathring K_\alpha$ is indexed in the standard way as $\{\mathring k_{\alpha,s} : s \in 2^{< \omega}\}$''.
		\item 
			$\one \forces_{\alpha+1}$ ``$\mathring g_\alpha = \bigcup \{ k_{\alpha, \mathring c _ \alpha \res n } : n \in \omega \}$''.
\end{enumerate}

Ignoring the ``$\Fn(\omega,2)$'', Conditions (1)(2)(3) are the standard setup for forcing $\MA$. We apply the usual bookkeeping to make sure that the $\mathring{\QQQ}_\alpha$ run through names for all possible c.c.c.\@ orders of size $< \kappa$; then ${\V}_\kappa$ satisfies $\MA + 2^{\aleph_0} = \kappa$. This all still works if we include the ``$\Fn(\omega,2)$'', which we use to construct the Gregory tree.

Conditions $(1-10)$ give us the Gregory tree $\T$ in ${\V}_\kappa$, named by $\mathring \T _\kappa$.  Condition (10) implies that $\T$ has no uncountable chains, and the usual bookkeeping would let us choose the $\mathring g _\alpha$ so that every Cantor subtree of $\T$ has a path. The main difficulty in the construction is in preserving (10). There are problems both at successors and at limits, addressed by Conditions $(11 - 15)$.

At successors:
Since $\mathring \T_{\alpha+1}$ is forced to be a subtree of $2^{< \omega_1}$, (9) requires $\mathring g_\alpha $ to be forced to be in $2^{< \omega_1}$, with all proper initial segments
of $\mathring g_\alpha $ in $\mathring \T_\alpha$. Since $\mathring \T_{\alpha+1} \backslash \mathring \T_\alpha$ is forced to be countable, Condition (10) is preserved in passing from
$\mathring \T_\alpha$ to $\mathring \T_{\alpha+1}$ \emph{unless} $\mathring \QQQ_\alpha$ adds a path through $\mathring \T_\alpha$, but this cannot happen by (11).

There is no problem ensuring (11) in the inductive construction.
$\Fn(\omega, 2)$ can never add a path through $\T_\alpha$.
To make sure that $\mathring \QQQ_\alpha$ does not add such a path, let $\mathring \QQQ_\alpha \cong \mathring \SSS_\alpha  * \mathring \RRR_\alpha$,
where $\mathring \SSS_\alpha$ is the name for the poset which specializes $\mathring \T_\alpha$.
Note that this does not interfere with the usual bookkeeping for making MA true.  Say this bookkeeping tells us that $\mathring \QQQ_\alpha$ should be $\mathring \BBB _\alpha$, which we may assume is always a $\PPP _\alpha$--name and that $\one \forces_\alpha$ ``$\mathring \BBB _\alpha$ is c.c.c.''; the c.c.c.\@ is not affected by the $\Fn(\omega,2)$, but it could be affected by the specializing order. Then $\mathring \RRR_\alpha$ is a $\PPP_\alpha * \Fn(\omega,2) * \mathring \SSS_\alpha$ name for the partial order which is $\mathring \BBB _\alpha$ if $\mathring \BBB _\alpha$ remains c.c.c.\@ after forcing with $\mathring \SSS_\alpha$, and otherwise
is the trivial order $\{\one\}$.

At limits:
In (8), we are literally taking the union of names in the ground model. This clearly preserves (4)(5) for $\T_\gamma$, and (10) is also preserved \emph{unless} $\cf(\gamma) = \omega_1$, in which case (10) might fail. For example, the $g_\alpha$ for $\alpha < \omega_1$ might all be compatible, yielding an uncountable chain in $\T_{\omega_1}$.

We avoid this problem by $(12 - 15)$. These say that working in ${\V}_{\alpha+1}$,
we choose the node $g_\alpha \in \T_{\alpha +1}$ as follows: We take a Cantor tree $K_\alpha \subseteq \T_\alpha$ (given to us by the usual bookkeeping) and let $g_\alpha$ be the path through this Cantor tree indexed by the Cohen real $c_\alpha$ added into ${\V}_{\alpha+1}$ by the $\Fn(\omega,2)$. Since $K_\alpha \in {\V}_\alpha$ by (13),
	\begin{equation*}
		g_\alpha \notin {\V}_\alpha. \tag{a}
	\end{equation*}
Now, suppose that (10) should fail at some point during the construction. Then we have $\gamma \le \kappa$ such that (10) holds for all $\alpha < \gamma$ but (10) fails at $\gamma$, so that we have a $\PPP_\gamma$--name $\mathring h$ and a $p \in \PPP_\gamma$ which forces that 
$\mathring h \in 2^{\omega_1}$ and is a path through $\mathring \T_\gamma$;
we may assume that $\mathring h$ is a nice name for a
subset of $\omega_1 \times 2$.
As noted above, $\gamma$ is a limit of cofinality $\omega_1$. Now, we argue both in $\V$ and in ${\V}[G]$, where $p \in G$ and $G$ is $(\V, {\P}_{\gamma})$--generic.

In $\V$, let $\langle \alpha_\xi : \xi < \omega_1 \rangle$ be a continuously increasing sequences of limit ordinals with supremum $\gamma$.  For $\mu < \omega_1$, we regard $\mathring h \res \mu$ as a nice $\PPP_\gamma$--name for an element of $2^\mu$; since $\PPP_\gamma$ is c.c.c., this $\mathring h \res \mu$ is actually a $\P_{\alpha_\xi}$--name for some $\xi < \omega_1$.  Then there is a club $C_0 \subset \omega_1$ such that $\mathring h \res \alpha_\xi$ is a $\PPP_{\alpha_\xi}$--name for each $\xi \in C_0$; so, in ${\V}[G]$, we have
	\begin{equation*} 
		h\res \alpha_\xi \in {\V}_{\alpha_\xi}. \tag{b}
	\end{equation*}
Also in $\V[G]$, each $\T_{\alpha_\xi}$  is special, so there is an $\eta > \xi$ such that 
$h \res \alpha_\eta \notin \T_{\alpha_\xi}$. Since we are taking unions of the trees at limit ordinals, there is a club $C_1 \subset \omega_1$ such that for $\xi \in C_1$ we have
	\begin{equation*} 
		h\res \alpha_\xi \notin \T_{\alpha_\xi}. \tag{c}
	\end{equation*}
Fix a limit ordinal $\xi \in C_0 \cap C_1$. Since $h \res \alpha_\xi \in \T_\gamma$, we may fix $\delta$ with $\alpha_\xi \le \delta < \gamma$ such that  $h \res \alpha_\xi \in \T_{\delta+1} \backslash \T_\delta$, which implies, by (9), that $h\res \alpha_\xi = {g_\delta}^{\frown}s$ for some $s \in 2^{< \omega}$, so $g_\delta \in {\V}_\delta$ by (b), contradicting (a).
\end{proof}

\section{Consistency of no Gregory Trees with Large Continuum}
\label{sec:noGtreesplusclarge}
In this section we shall prove:

\begin{Theorem}
\label{Thm:isnt-greg-tree}
Assume that in the ground model $\V$:
\begin{itemize}
\item[1.]
$\kappa \ge \aleph_2$ and $\kappa^{<\kappa} = \kappa$.
\item[2.]
$ \lambda^{\aleph_0} < \kappa$ for all $\lambda < \kappa$.
\item[3.]
${\diamondsuit}_{\kappa}(S)$, where
$S = \{\alpha < \kappa: \cf(\alpha) = \omega_1\}$.
\end{itemize}
Then there is a c.c.c.\@ forcing
extension $\V[G]$ satisfying $\MA + 2^{\aleph_0} = \kappa$ in which
there are no Gregory trees. 
\end{Theorem}

We do not know whether (3) follows from (1)(2); it does
by Gregory \cite{Gr2} in the case that $\kappa = \lambda^+$
and $\lambda^{\aleph_1} = \lambda$.
If we start with $\V = \LL$, then (1) and (3) hold for all regular
$\kappa \ge \aleph_2$, but (2) fails if $\kappa$ is the successor
to a cardinal of cofinality $\omega$, so we are left with 
Question \ref{ques:greg}.

As with the proof of Theorem \ref{Thm:is-greg-tree},
we shall modify the usual ccc iteration to produce a model
of $\MA + 2^{\aleph_0} = \kappa$ (using (1)).
To kill a potential Gregory tree $\T$ in $\VG$, we use (2) plus
countably closed elementary submodels to produce a club $C\subseteq \kappa$
such that $\T ^ \alpha := T \cap \V[G_\alpha]$
has the CTP in $\V[G_\alpha]$ for all $\alpha \in C\cap S$.
Then, we use (3) to ensure that at some stage $\alpha$
in the construction, we kill $\T ^ \alpha$ by shooting
a cofinal branch through it, so that we also kill $\T$.

Now, to kill $\T ^ \alpha$ by a c.c.c.\@ poset, we cannot force
with $\T ^ \alpha$, since this is not c.c.c.
Instead, we shall find a Suslin subtree $\Q_\alpha \subset \T ^ \alpha$
and force with $\Q_\alpha$.
This method is patterned after \cite{KH}, which proved 
Theorem \ref{Thm:isnt-greg-tree} in the special case
that $\kappa = \aleph_2$ and $\diamondsuit$ (that is, $\diamondsuit_{\omega_1}$)
holds in $\V$.  It is well-known that $\diamondsuit$ will
remain true in $\V[G_\alpha]$ (since $\alpha < \omega_2$),
and hence, by Lemma 5.8 of \cite{KH}, the tree
$\T ^ \alpha$ will have a Suslin subtree.
But for longer iterations, 
$\diamondsuit$ (and CH) will fail whenever $\alpha \geq \omega_2$.
Instead, we shall use the fact that $\cf(\alpha) = \omega_1$.
It is well-known that this implies that there is a Suslin tree
in $\V[G_\alpha]$, since Cohen reals have been added cofinally often below 
$\alpha$ (see, for example, Theorems 3.1 and 6.1 of \cite{MHD}).
Here, we shall prove 
Theorem \ref{thm:Suslinsubtreesatcfomega1}, which shows how to
get the Suslin tree $\Q_\alpha$ inside of $\T ^ \alpha$.

\begin{Def}
For a limit ordinal $\gamma$ and a subtree $\T \subseteq 2^{< \gamma}$:
$\T$ is \emph{uniformly of height} $\gamma$ iff
$\forall f \in \T \, \forall \alpha < \gamma \,
\exists h \in \T \, [h < f \ \&\ \hgt(h) > \alpha ]$, and
$\T$ is \emph{branchy} iff
$f \cat 0, f \cat 1 \in \T$ for all $f \in \T$.  If $\T \in \V$
and $g : \gamma \to 2$, then $g$ is \emph{$\T$--generic over} $\V$
iff $\{g \res \alpha : \alpha < \gamma\}$ is $\T$--generic over $\V$.
\end{Def}

Such a tree is an atomless forcing order, and every
$\T$--generic filter is a path through $\T$.  If $\T$
is countable, then $\T$ is equivalent to Cohen forcing
$\Fn(\omega,2)$. 
We can now modify the standard Jensen construction of
a Suslin tree $\T \subseteq 2^{< \omega_1}$;
the Cohen reals allow us to replace the
use of $\diamondsuit$ at limits $\gamma < \omega$ by the requirement
that all $g \in \T \cap 2^\gamma$ be $\T \cap 2^{< \gamma}$--generic.
This is described in Lemma \ref{lemma:get-suslin}, which we shall prove
after listing some further conventions for names in c.c.c.\@ forcing extensions.

Say $\P$ is c.c.c.\@ and
$p \forces \mathring a \in 2^{< \omega_1}$.  
Then $p$ may not decide what the height (= domain)
$\hgt(\mathring a)$ is, but there is a $\xi < \omega_1$ such
that $p \forces \hgt(\mathring a) \le \xi$, so $\mathring a$
is forced to be a subset of $\xi \times 2$, and there
is a nice name  $\mathring b$ for a subset of $\xi \times 2$ such
that $p \forces \mathring a  = \mathring b$.

Next, consider subsets $A \subseteq 2^{< \omega_1}$ in $\VG$;
$A$ may be a tree, or an antichain in a tree; again, $\P$ is c.c.c.
$A$ is not a subset of a ground model set, 
but we may simplify the name for $A$ as follows.  Say $p$ forces that
$\mathring A \subseteq 2^{< \omega_1}$
and $1 \le |\mathring A| \le \kappa$.  Then, in $\VG$, we may list
$A$ in a $\kappa$--sequence (possibly with repetitions),
so there is a name $\mathring B$ such that
$p \forces \mathring A  = \mathring B$, where
$\mathring B = \{ \langle \mathring b _ \mu, p \rangle : \mu  < \kappa \}$
and each $\mathring b _ \mu$ is a nice name for a subset of
some $\xi_\mu \times 2$, where $\xi_\mu < \omega_1$ and
$\one \forces \mathring b _ \mu \in  2^{< \omega_1} \wedge
\hgt(\mathring b _ \mu ) \le \xi_\mu$.

These conventions make it easy to apply elementary submodel
arguments in $\V$ to the forcing construction.  For example,

\begin{Lemma}
\label{lemma:reflect-ctp}
Assume in $\V$: 
$\kappa \ge \aleph_2$ is regular and
$ \lambda^{\aleph_0} < \kappa$ for all $\lambda < \kappa$
and
$\langle \PPP_\alpha : \alpha \le \kappa\rangle$ is
a normal chain of c.c.c.\@  posets,
with each $|\P_\alpha| < \kappa$.  Let $\mathring \T =
\{ \langle \mathring b _ \mu, \one \rangle : \mu  < \kappa \}$
be a $\P_\kappa$--name which is forced by
$\one$ to be a subtree of $2^{< \omega_1}$ with the CTP, 
where each $\mathring b _ \mu$ is a nice name for a subset of
some $\xi_\mu \times 2$ with $\xi_\mu < \omega_1$.
Let $\mathring \T^\alpha =
\{ \langle \mathring b _ \mu, \one \rangle : \mu  < \alpha \}$.

There is then a club $C \subseteq \kappa$ such that for all $\alpha \in C$ with
$\cf(\alpha) > \omega$: $\mathring \T^\alpha$ is a $\P_\alpha$--name
and  $\one \forces_\alpha$ ``$\,\mathring \T^\alpha$ is
a subtree of $2^{< \omega_1}$ with the CTP''. 
\end{Lemma}
\begin{proof}
Fix a suitably large regular $\theta$.
Given the assumptions on $\kappa$, it is sufficient to prove
that the conclusion to the lemma holds whenever
$\alpha$ is an ordinal of the form $M \cap \kappa$,
where $M \prec H(\theta)$ is a countably
closed elementary submodel containing the relevant objects.

The fact that $\mathring \T^\alpha$ is a $\P_\alpha$--name is immediate
and does not need countable closure.  Likewise, to show that
$\one \forces_\alpha$ ``$\,\mathring \T^\alpha$ is
a subtree'', note that for each $\mu < \kappa$ there is a countable
$R \subseteq \kappa$ such that 
$\one\forces \forall \eta \le \xi_\mu \exists \nu \in R
\, [\mathring b_\mu \res \eta = \mathring b_\nu]$, and
by $M \prec H(\theta)$, some such $R$ is in $M$, so that $R \subset \alpha$.
The proof of the CTP is similar, but uses the countable closure
of $M$ to imply that $M$ contains $\P_\alpha$--names for
every possible Cantor subtree of $\mathring \T^\alpha$
which lies in $\V_\alpha$.
\end{proof}

Similar (and easier) reflection arguments work for sets $A$
of size $\aleph_1$.
Call $A \subseteq 2^{< \omega_1}$  \emph{skinny} iff $|A| = \aleph_1$ and
each $A \cap 2^\xi$ is countable.  Then we can list $A$
in an $\omega_1$--sequence, listing nodes in order of their height.
If $p$ forces that $\mathring A$ is a skinny subset of 
$\omega_1$, then there is a club $C$ and a name $\mathring B$ such that
$p \forces \mathring A  = \mathring B$, where $\mathring B =
\{ \langle \mathring b _ \mu, p \rangle : \mu  < \omega_1 \}$ as above,
and also $\one \forces \hgt(\mathring b _ \mu ) \ge \gamma$ 
whenever $\mu \ge \gamma \in C$.
If $\mathring B _\gamma$ is the name
$\{ \langle \mathring b _ \mu, p \rangle : \mu  < \gamma \}$,
then $\one \forces \mathring B \cap 2^{< \gamma} = \mathring B _\gamma$
whenever $\gamma \in C$.
With iterated forcing, where $\P = \P_{\omega_1}$ results from a
normal chain of c.c.c.\@ posets
$\langle \PPP_\alpha : \alpha \le \omega_1 \rangle$,
we can also arrange for $ \mathring B _\gamma$ to be 
a $\P_\gamma$--name whenever $\gamma \in C$, so that from the 
point of view of $\VG$ with $p \in G$, each
$A \cap 2^{< \gamma} \in \V[G_\gamma]$.

\begin{Lemma}
\label{lemma:get-suslin}
Suppose that in $\V$, 
$\langle \PPP_\alpha : \alpha \le \omega_1 \rangle$ is
a normal chain of c.c.c.\@ posets and $G$ is $\PPP_{\omega_1}$--generic.
In $\VG$, suppose that the subtree $\T \subseteq 2^{< \omega_1}$
is uncountable, branchy, skinny, and uniformly of height $\omega_1$.
Assume also that there is club of limit ordinals
$C \subset \omega_1$ such that for all $\gamma \in C$:
$\T \cap 2^{< \gamma} \in \V[G_\gamma]$ and every 
$g \in \T \cap 2^\gamma$ is $\T \cap 2^{< \gamma}$--generic
over $\V[G_\gamma]$.  Then $\T$ is Suslin in $\VG$.
\end{Lemma}
\begin{proof}
If not, then in $\VG$ we have an uncountable maximal antichain $A \subseteq \T$.
Then there is a club $C_1$ such that $A \cap 2^{< \gamma}$ is 
a maximal antichain in $\T \cap 2^{< \gamma}$ for all $\gamma \in C_1$.
Also, $A$ is skinny, so as pointed out above, there is a club
$C_2$ such that $A \cap 2^{< \gamma} \in \V[G_\gamma]$ for
all $\gamma \in C_2$.  Now, fix $\gamma \in C \cap C_1 \cap C_2$,
and consider any $g \in \T \cap 2^\gamma$.  Then
$\{g \res \xi : \xi < \gamma\}$ is generic over $\V[G_\gamma]$,
so it meets the maximal antichain $A \cap 2^{< \gamma} \in \V[G_\gamma]$.
But then $A \subseteq 2^{< \gamma}$, so $A$ is countable.
\end{proof}

This lemma lets us construct a Suslin tree in an iterated
forcing extension, but we actually need our tree to be inside
a given Gregory tree.  To do that, we use the following lemma,
which is related to the well-known fact that adding a Cohen real
actually adds a perfect set of Cohen reals:

\begin{Lemma}
\label{lemma:get-perfect}
Assume that in $\V$:  $\gamma$ is a countable limit ordinal
and $\T$ is a countable subtree of 
$2^{< \gamma}$ which is branchy and uniformly of height $\gamma$.
Fix an $\omega$--sequence $\langle \alpha_i : i < \omega \rangle$
of ordinals increasing to $\gamma$ and fix any $f\in \T$ with
$\hgt(f) = \alpha_0$.
Let $\VG$ be \emph{any} forcing extension of $\V$ which contains a Cohen real
\textup(e.g., a filter which is $\Fn(\omega,2)$--generic over $\V$\textup).

Then, in $\VG$:  There is a Cantor tree
$\{f_{\sigma}: \sigma \in 2^{< \omega}\} \subseteq \T$ such
that $f_\emptyset = f$ and $\alpha_i \le \hgt(f_\sigma) < \gamma$ whenever
$\hgt(\sigma) = i$, and such that for all $\psi \in 2^\omega$,
$\bigcup\{f_{\psi\res i} : i < \omega\}$ is
$\T$--generic over $\V$.
\end{Lemma}
\begin{proof}
In $\V$, let $\P$ be the poset of
``partial Cantor trees starting at $f$''.
So, $p \in \P$ iff for some $n = n_p \in \omega$:
$p$ is a function from
$ 2^{\leq n}$ into $\T$, $p(0) = f$,
 $\alpha_i \le \hgt(f_\sigma) < \gamma$
whenever $\hgt(\sigma) = i \le n$, and  
$p({\sigma}^{\frown}{0}) \perp p({\sigma}^{\frown}{1})$   
whenever $\hgt(\sigma) < n$.
Order $\P$ by $q \le p$ iff $q \supseteq p$.
$\P$ is countable, so the existence of a Cohen real implies the
existence of a $(\V,\P)$--generic filter $H \in \VG$.
Note that $\{p \in \P: n_p \ge m\}$
is dense for each $m$ because
$\T$ is uniformly of height $\gamma$,
so $\bigcup H : 2^{<\omega} \to \T$.
Let $f_\sigma = (\bigcup H )(\sigma)$.
Fix any $\psi \in 2^\omega$.  To verify
$\T$--genericity of $\bigcup\{f_{\psi\res i} : i < \omega\}$,
let $D \subseteq \T$ be dense.  Then
$D^* := \{p: \forall \sigma \in 2^{n_p}\, [p(\sigma) \in D]\}$
is dense in $\P$.
If $p \in H \cap D^*$ then $f_{\psi \res n_p} \in D$.
\end{proof}

In any non-trivial iterated forcing, the Cohen reals come for free because
of the following well-known lemma.

\begin{Lemma}
\label{lemma:get-Cohen}
Suppose that in $\V$,  $\gamma$ is any limit ordinal and
$\langle \PPP_\alpha : \alpha \le \gamma \rangle$ is
a normal chain of c.c.c.\@ posets.
Let $G$ be $\PPP_\gamma$--generic over $\V$, and assume that
$\VG \ne \V[G_\alpha]$ for any $\alpha < \gamma$.
Then $\VG$ contains a real which is Cohen generic over $\V$.
\end{Lemma}

This lemma is actually not critical for our proof, since in iterating
to make $\MA$ true, we could easily add a Cohen real explicitly at
each stage.

\begin{Theorem}
\label{thm:Suslinsubtreesatcfomega1}
Suppose that in $\V$:  $\pi$ is a limit ordinal with
$\cf(\pi) = \omega_1$, and
$\langle \PPP_\alpha : \alpha \le \pi \rangle$ is
a normal chain of c.c.c.\@ posets.

Let $G$ be $\PPP_\pi$--generic over $\V$, and assume that
$\VG \ne \V[G_\alpha]$ for any $\alpha < \pi$.

In $\VG$:  Let $\T$ be a subtree of $2^{< \omega_1}$ with the CTP.
Then $\T$ has a Suslin subtree.
\end{Theorem}
\begin{proof}
First, restricting to a club and applying Lemma \ref{lemma:get-Cohen}
(using the various $\V[G_\alpha]$ as the ground model),
we may assume that $\pi = \omega_1$ and that 
each $\V[G_{\alpha+1}]$ contains a real which is Cohen generic
over $\V[G_\alpha]$.  Actually, the club is obtained in $\V[G]$;
but it then contains a ground model club by the c.c.c.; so the
restricted sequence of forcing posets also lies in $\V$.

Now, working in $\V[G]$, we construct a Suslin subtree $\SSS \subset \T$
by constructing inductively $\SSS \cap 2^{<\gamma}$.
$\SSS$ will be branchy, so we only need to specify the construction
for limit $\gamma$.  We assume that we have $\SSS \cap 2^{<\gamma}$,
and we assume (inductively) that each such $\SSS \cap 2^{<\gamma}$
is countable and is uniformly of height $\gamma$, and we must describe
$\SSS \cap 2^{\gamma}$.
For each $f \in \SSS \cap 2^{<\gamma}$, choose a
$g_f \in 2^{\gamma}$ such that $g_f < f$ and such that
$g_f \res \xi \in \SSS \cap 2^{<\gamma}$ for all $\xi < \gamma$;
then $\SSS \cap 2^{\gamma} = \{g_f : f \in \SSS \cap 2^{<\gamma}\}$.
To get each $g_f$:
Fix an $\omega$--sequence $\langle \alpha_i : i < \omega \rangle$
of ordinals increasing to $\gamma$, with $\alpha_0 = \hgt(f)$.
Then choose a Cantor tree
$\{f_{\sigma}: \sigma \in 2^{< \omega}\} \subseteq \SSS \cap 2^{<\gamma}$ such
that $f_\emptyset = f$ and
$\alpha_i \le \hgt(f_\sigma) < \gamma$ whenever $\hgt(\sigma) = i$;
this is easily done since $\SSS$ is uniformly of height $\gamma$.
Furthermore, \emph{if} $\SSS \cap 2^{<\gamma} \in \V[G_\gamma]$,
apply Lemma \ref{lemma:get-perfect} and assume that
for all $\psi \in 2^\omega$,
$\bigcup\{f_{\psi\res i} : i < \omega\}$ is
$(\SSS \cap 2^{<\gamma})$--generic over $\V[G_\gamma]$.
Whether or not $\SSS \cap 2^{<\gamma} \in \V[G_\gamma]$, we apply the CTP
of $\T$ to always choose $g_f = \bigcup_i f_{\psi\res i} \in \T$.

Now that we have constructed $\SSS \subset \T$
in $\VG$, we note that it is skinny,
so as pointed out above, there is a club of limits $C \subset \omega_1$
such that $\SSS \cap 2^{<\gamma} \in \V[G_\gamma]$ for all $\gamma \in C$.
For these $\gamma$,  every
$g \in \SSS \cap 2^\gamma$ is $\SSS \cap 2^{< \gamma}$--generic
over $\V[G_\gamma]$, so,
by Lemma \ref{lemma:get-suslin}, $\SSS$ is Suslin in $\VG$.
\end{proof}

\begin{proof}[Proof of Theorem \ref{Thm:isnt-greg-tree}]
As with the proof of Theorem \ref{Thm:is-greg-tree}, in $\V$ we
inductively construct the following, satisfying the listed conditions:
\begin{enumerate}
\item 
$\langle \PPP_\alpha : \alpha \le \kappa\rangle$ is a normal chain of c.c.c.\@ posets.
\item 
$|\PPP_\alpha| < \kappa$ for all $\alpha < \kappa$.
\item 
$\PPP_{\alpha + 1} \cong \PPP_\alpha *  \mathring{\QQQ}_\alpha$,
where $\one \forces_\alpha$ ``$\mathring{\QQQ}_\alpha$ is c.c.c.''.
\item 
$\PPP_\alpha$ as a \emph{set} is the ordinal $\delta_\alpha\le \kappa$,
with $\one = 0$.
\end{enumerate}
As before, (1)(2)(3) are the standard setup for forcing $\MA$.
We apply the usual bookkeeping to make sure that the
$\mathring{\QQQ}_\alpha$, for $\cf(\alpha) \ne \omega_1$,
run through names for all possible atomless c.c.c.\@ orders of size $< \kappa$,
so $\VG = {\V}_\kappa$ satisfies $\MA + 2^{\aleph_0} = \kappa$. 
Condition (4) is irrelevant for this, although it is sometimes included
in expositions to facilitate the bookkeeping.
Note that $\kappa$ is regular, so the $\delta_\alpha$, for $\alpha < \kappa$,
form a continuously increasing sequences of ordinals less than $\kappa$,
and $\delta_\kappa = \kappa$;
also note that $\{\alpha < \kappa : \delta_\alpha = \alpha\}$ is a club.
We have included (4) to facilitate the use of 
${\diamondsuit}_{\kappa}(S)$, which will give us
$\Q_\alpha$ when $\cf(\alpha) = \omega_1$.

To show that there are no Gregory trees in $\VG$,
it is sufficient to show in $\V$ that whenever
$\one$ forces $\mathring \T$ to be a subtree of $2^{< \omega_1}$ with the CTP,
$\one$ also forces $\mathring \T$ to have a cofinal branch.
By the CTP, $|\T| =  2^{\aleph_0} = \kappa$ in $\VG$, so as noted above, 
$\T$ has a name of the form
$\mathring \T = \{ \langle \mathring b _ \mu, \one \rangle : \mu  < \kappa \}$,
where each $\mathring b _ \mu$ is a nice name for a subset of
some $\xi_\mu \times 2$, where $\xi_\mu < \omega_1$ and
$\one \forces \mathring b _ \mu \in  2^{< \omega_1} \wedge
\hgt(\mathring b _ \mu ) \le \xi_\mu$.

We must specify our $\diamondsuit$ sequence
before we have defined an order on the sets $\P_\alpha = \delta_\alpha$.
The definition of $\check x$ only uses the identity $\one = 0$,
but the notion of ``nice name'' presupposes that we know what
an antichain is.  So,
call $\mathring b$ a \emph{pseudo-nice} $\delta$--name for a subset of 
$X \in V$ iff
$\mathring b = \bigcup\{ \{\check x\} \times E_x : x \in X \}$, where
each $E_x \in [\delta]^{\le \omega}$.
Then every nice name using the eventual
order on $\delta_\alpha$ will be also pseudo-nice.

Our  $\diamondsuit$ sequence will make believe that $\delta_\alpha = \alpha$,
since this is true on a club.
So, for $\alpha \in S  = \{\alpha < \kappa: \cf(\alpha) = \omega_1\}$,
choose a $\mathring \T_\alpha$ of the form
$ \{ \langle \mathring b^\alpha _ \mu, \one \rangle : \mu  < \alpha \}$,
where each $\mathring b^\alpha _ \mu$ is a pseudo-nice $\alpha$--name for a
subset of some $\xi_\mu^\alpha \times 2$, where $\xi_\mu^\alpha < \omega_1$.
These $\mathring \T_\alpha$ must have the $\diamondsuit$ property that
whenever $\mathring \T =
\{ \langle \mathring b _ \mu, \one \rangle : \mu  < \kappa \}$
has the analogous form (replacing $\alpha$ with $\kappa$),
the set of $\alpha \in S$ for which
$\mathring \T_\alpha =
\{ \langle \mathring b _ \mu, \one \rangle : \mu  < \alpha \}$
is stationary.

Now, when $\alpha \in S$ and we have constructed $\P_\alpha$
(i.e., we know the ordinal $\delta_\alpha$ and its ordering),
choose $\mathring \Q _\alpha$ as follows:
$\mathring \Q _\alpha$ is a name for the trivial one-element order 
\emph{unless} $\delta_\alpha = \alpha$ and 
each $\mathring b _ \mu^\alpha$ is indeed a nice $\P_\alpha$--name and
$\one \forces_\alpha \mathring b _ \mu^\alpha \in  2^{< \omega_1}$
and $\one \forces_\alpha$ ``$\mathring \T_\alpha$ is 
a subtree of $2^{< \omega_1}$ with the CTP''.  In that case,
$\V_\alpha$ will contain the tree $\T_\alpha$ which (in $\V_\alpha$)
has the CTP, and then Theorem \ref{thm:Suslinsubtreesatcfomega1} applies
to construct a Suslin subtree $\Q_\alpha \subset \T_\alpha$.  Then,
back in $\V$, we let $\mathring \Q _\alpha$ be a name for this
$\Q_\alpha$, so that in $\V_{\alpha+1}$ we have a cofinal branch in $\Q_\alpha$.

Finally to show that there are no Gregory trees in $\VG$,
assume that $\one$ forces $\mathring \T$ to be a subtree
of $2^{< \omega_1}$ with the CTP. 
Let the $\mathring \T^\alpha$ be as in Lemma  \ref{lemma:reflect-ctp}.
Then, by Lemma  \ref{lemma:reflect-ctp},
there is then a club $C$ of $\omega_1$--limits such that for $\alpha \in C$:
$\delta_\alpha = \alpha$,
$\mathring \T^\alpha$ is a $\P_\alpha$--name,
and  $\one \forces_\alpha$ ``$\,\mathring \T^\alpha$ is
a subtree of $2^{< \omega_1}$ with the CTP''. 
Choosing $\alpha \in C$ with $\mathring \T^\alpha = \mathring \T_\alpha$
shows that $\one$ forces that there is (in $V_{\alpha + 1})$
a cofinal branch in $\mathring \T $.
\end{proof}

\end{document}